\documentclass[11pt]{amsart}
\usepackage[T1]{fontenc}
\usepackage{a4wide}
\usepackage[latin1]{inputenc}
\usepackage{amsmath}
\usepackage{amsfonts}
\usepackage{amssymb}
\usepackage{amsthm}
\usepackage{subfigure}
\newtheorem{theo}{Theorem}[section]
\newtheorem{prop}[theo]{Proposition}
\newtheorem{coro}[theo]{Corollary}
\newtheorem{lemm}[theo]{Lemma}

\theoremstyle{definition}
\newtheorem{def1}[theo]{Definition}
\theoremstyle{remark}
\newtheorem*{rema}{Remark}
\newcommand{\Op}{{\rm{Op}}}
\newcommand{\Vol}{{\rm{Vol}}}
\newcommand{\R}{\mathbb{R}} 
\newcommand{\Z}{\mathbb{Z}} 
\newcommand{\N}{\mathbb{N}} 
\newcommand{\cD}{\mathcal{D}}
\newcommand{\cC}{\mathcal{C}}
\newcommand{\cO}{\mathcal{O}}
\newcommand{\cP}{\mathcal{P}}
\newcommand{\cS}{\mathcal{S}}

\newcommand{\cL}{\mathcal{L}}
\newcommand{\To}{\longrightarrow} 
 \newcommand{\Lim}{\mathop{\longrightarrow}\limits}
 \newcommand{\Tr}{{\rm Tr}}
 \newcommand{\bbbone}{{\mathchoice {1\mskip-4mu {\rm{l}}} {1\mskip-4mu {\rm{l}}}
{ 1\mskip-4.5mu {\rm{l}}} { 1\mskip-5mu {\rm{l}}}}}

\title[Schr\"odinger equation on negatively curved manifolds]{Dispersion and controllability for the Schr\"odinger equation on negatively curved manifolds}
\author{Nalini Anantharaman}
\address{Universit\'e Paris Sud, Laboratoire de Math\'ematique (UMR 8628), B\^atiment 425, 91405 Orsay Cedex, France}
\email{Nalini.Anantharaman@math.u-psud.fr}

\author{Gabriel Rivière}
\address{\'Ecole Polytechnique, Centre de Math\'ematiques Laurent Schwartz (UMR 7640), 91128 Palaiseau Cedex, France}
\email{gabriel.riviere@polytechnique.edu}

\thanks{N. Anantharaman wishes to acknowledge the support of Agence Nationale de la Recherche,
under the grant ANR-09-JCJC-0099-01.}

\begin{document}

\begin{abstract}
We study the time-dependent Schr\"odinger equation $
\imath\frac{\partial u}{\partial t}=-\frac{1}{2}\Delta u,$ on a compact riemannian manifold on which the geodesic flow has the Anosov property.
Using the notion of semiclassical measures, we prove various results related to the dispersive properties of the Schr\"odinger propagator, and to controllability.
\end{abstract}
\maketitle

\section{Introduction}

Let $M$ be a smooth compact riemannian manifold of dimension $d$ (without boundary). 
 We denote by $\Delta$ the laplacian on $M$. We are interested in understanding the regularizing
properties of the Schr\"{o}dinger equation
\[
\imath\frac{\partial u}{\partial t}=-\frac{1}{2}\Delta u,\qquad u\rceil_{t=0}\in
L^{2}(M).
\]
More precisely, given a sequence of initial conditions $u_{n}\in
L^{2}(M)$, we investigate the asymptotic behaviour of the family of probability measures
\begin{equation}
\nu_{n}(dx)=\left(  \int_{0}^T|e^{\imath t\Delta/2}u_{n}(x)|^{2}dt\right)  d{\rm Vol}(x)
\label{probmes}%
\end{equation}
(where ${\rm Vol}$ denotes the riemannian volume measure on $M$).

We want to relate this question to the behaviour of the geodesic flow, using results on propagation of singularities. For that purpose, we reformulate the question using the notion of semiclassical measures.
We consider a sequence of states $(u_{\hbar})_{\hbar\rightarrow0^+}$ normalized in $L^2(M)$ (indexed by a parameter $\hbar>0$ going to $0$, which plays the role of Planck's constant in quantum mechanics), and for every $t\in\R$ we define the following family of distributions on the cotangent bundle $T^*M$:
\begin{equation}\label{general-def-distrib}
\forall a\in\mathcal{C}^{\infty}_o(T^*M), {\mu}_{\hbar}(t)(a)=\int_{T^*M}a(x,\xi)d{\mu}_{\hbar}(x,\xi):=\langle e^{\imath t\Delta
/2} u_{\hbar}|\Op_{\hbar}(a)| e^{\imath t\Delta
/2}u_{\hbar}\rangle_{L^2(M)},
\end{equation}
where $\Op_{\hbar}(a)$ is a $\hbar$-pseudodifferential operator of principal symbol $a$ (see ~\cite{DS}, or appendix~\ref{pdo-manifold} for a brief reminder). This construction gives a description of a state in terms of position and impulsion variables.  
Throughout the paper, we will denote $U^t:=e^{\imath t\Delta/2}$ the quantum propagator.

By standard estimates on the norm of $\Op_{\hbar}(a)$ (the
Calder\'{o}n-Vaillancourt theorem), the map $t\mapsto
{\mu}_{\hbar}(t)$ belongs to $L^{\infty}(\R;\cD^{\prime}\left(  T^{\ast
}M\right)  )$, and is uniformly bounded in that space as
$\hbar\To0^{+}$. Thus, one can extract subsequences that converge in the
weak-$\ast$ topology on $L^{\infty}(\R;\cD^{\prime}\left(  T^{\ast}%
M\right)  )$. In other words, after possibly extracting a
subsequence, we have
\[
\mu_{\hbar}(a\otimes\theta)=\int_{\R}\theta(t)a(x,\xi){\mu}_{\hbar}(t)(dx,d\xi)dt\Lim_{\hbar\To0}\int_{\R}%
\theta(t)a(x,\xi)\mu(t)(dx,d\xi)dt
\]
for all $\theta\in L^{1}(\R)$ and $a\in \cC_{o}^{\infty}(T^{\ast}M).$
The main example to keep in mind is the case when $\theta$ is the characteristic function of some interval $[0, T]$. In that case we can write
$$\mu_{\hbar}(a\otimes\theta)=\int_0^T  \langle e^{\imath t\Delta
/2} u_{\hbar}|\Op_{\hbar}(a)| e^{\imath t\Delta
/2}u_{\hbar}\rangle dt=\hbar \int_0^{T/\hbar}  \langle e^{\imath t\hbar\Delta
/2} u_{\hbar}|\Op_{\hbar}(a)| e^{\imath t\hbar\Delta
/2}u_{\hbar}\rangle dt.$$
In the last term we have just expressed everything in terms of the flow $e^{\imath t\hbar\Delta
/2}$, which solves the equation 
$
-\frac{\hbar^2}{2}\Delta v =\imath\hbar\frac{\partial v}{\partial t}$ with the time-parametrization of quantum mechanics. Thus, in the time-scale of quantum mechanics, we are averaging
over time intervals of order $\hbar^{-1}$.

It follows from standard properties of $\Op_\hbar(a)$ that
the limit $\mu$ has the following properties~:

\begin{itemize}
\item 
for almost all $t$, $\mu(t)$ is a positive measure on $T^{\ast
}M$.

\item the unitarity of $U^t$ implies that $\int_{T^{\ast}%
M}\mu(t)(dx,d\xi)$ does not depend on $t$; from the normalization
of $u_{\hbar}$, we have $\int_{T^{\ast}M}\mu(t)(dx,d\xi)\leq1$, the
inequality coming from the fact that $T^{\ast}M$ is not compact,
and that there may be an escape of mass to infinity.

\item define the geodesic flow $g^{s}:T^{\ast}M\To T^{\ast
}M$ as the hamiltonian flow associated with the energy $p(x, \xi)=\frac{||\xi||_x^2}2$. From the \emph{Egorov theorem}, we have
\[
\forall s\in\mathbb{R},\ e^{-is\hbar\Delta/2}\Op_{\hbar}(a)e^{is\hbar\Delta/2}=\Op_{h}(a\circ g^{s})+O_{s, a}(\hbar)
\]
for $a\in C_o^\infty(T^*M)$.
At the limit $\hbar\To0^{+}$, this implies that $\mu(t)$ is invariant under $g^{s}$, for almost all
$t$ and all $s$.
\end{itemize}

These sequences of distributions were already studied by Macià \cite{Mac2}; we refer to that paper for details about the facts mentioned above. Maci\`a was mostly interested in describing the properties of the measures $\mu(t)$ in the case where the geodesic flow on the manifold $M$ was not chaotic (Zoll manifolds for instance, or the flat torus \cite{Mac3, AM}).

In this paper, we are interested in a completely different situation where the geodesic flow has the Anosov property (manifolds of negative curvature are the main example). In this setting, the case where the initial states $u_\hbar$ are eigenfunctions of the laplacian, satisfying $-\hbar^2\Delta u_\hbar=u_\hbar$, has been much studied; in this particular situation $\mu_\hbar(t)$ does not depend on $t$. The Shnirelman theorem (also called Quantum Ergodicity Theorem) says that for ``most'' sequences of eigenfunctions
$u_\hbar$, the limit $\mu$ is the Liouville measure on the unit cotangent bundle $S^*M$; see \cite{Sc, Ze, CdV} for the precise statement.
It is also known, by the work of Anantharaman and Nonnenmacher, that for any sequence of eigenfunctions the limit $\mu$ has positive entropy \cite{A, AN2}.
The aim of this paper is twofold: \emph{extend the Shnirelman theorem to the setting of the time dependent equation} and \emph{prove lower bounds on the metric entropy of the measures $\mu(t)$}. We shall also show how these results apply to the controllability problem for the Schr\"odinger equation.

\section{Statement of the results}

\subsection{Semiclassical large deviations} Our first result is a generalization (and a reinforcement in the case of Anosov geodesic flows) of the quantum ergodicity theorem. Recall that the Shnirelman Theorem is a result on orthonormal bases of eigenfunctions of the laplacian. In order to state an analogue of it for solutions of the time dependent Schr\"odinger equation, we introduce a notion of generalized orthonormal families.


\subsubsection{Generalized orthonormal family} We fix $\alpha>0$ and a sequence $I(\hbar):=[a(\hbar),b(\hbar)]$ of subintervals that are of length at least $2\alpha\hbar$ for every $\hbar>0$. We also suppose that $\lim_{\hbar\rightarrow0^+}a(\hbar)=\lim_{\hbar\rightarrow0^+}b(\hbar)=1.$ We denote $N(I(\hbar))$ the number of eigenvalues $\lambda_j^2$ of $\Delta$ (counted with their multiplicities) satisfying $\hbar^2\lambda_j^2\in I(\hbar)$. We assume that
\begin{equation}\label{e:weyl}N(I(\hbar))=\frac{\Vol(M)}{(2\pi\hbar)^d}\Vol( B_d(0, 1))(b(\hbar)-a(\hbar))(1+o(1))\end{equation}
(where $\Vol(M)$ is the riemannian volume of $M$, and $\Vol( B_d(0, 1))$ is the volume of the unit ball in $\R^d$).
According to~\cite{DuisGui}, we know that the Weyl's law \eqref{e:weyl} still holds in the case where $b(\hbar)-a(\hbar)=2\alpha\hbar$ if we suppose that the set of closed geodesics is of zero Liouville measure on $S^*M$ (this is the case for Anosov geodesic flows).

We introduce the notion of generalized orthonormal family localized in the energy window $I(\hbar)$:
\begin{def1}\label{d:generalized-orth-basis} For $\hbar>0$, let $(\Omega_{\hbar},\mathbb{P}_{\hbar})$ be a probability space and $u_{\hbar}:\Omega_{\hbar}\rightarrow L^2(M)$ a measurable map. We say that $(u_{\hbar}(\omega))_{\omega\in(\Omega_{\hbar},\mathbb{P}_{\hbar})}$ is a generalized orthonormal family (G.O.F.) of the spectral window $I(\hbar)$ if
\begin{itemize}
\item $\|u_{\hbar}(\omega)\|_{L^2(M)}=1+o(1)$ as $\hbar$ tends to $0$ (uniformly for $\omega$ in $\Omega_{\hbar}$);
\item $\displaystyle\left\|\left(\text{Id}_{L^2(M)}-\bbbone_{I(\hbar)}\left(-\hbar^2\Delta\right)\right)u_{\hbar}(\omega)\right\|_{L^2(M)}=o(1)$ as $\hbar$ tends to $0$ (uniformly for $\omega$ in $\Omega_{\hbar}$);
\item for every $B$ in $\mathcal{L}(L^2(M))$,
\begin{equation}\label{e:generalized-partitions-M}
\int_{\Omega_{\hbar}}\langle u_{\hbar}(\omega)|B|u_{\hbar}(\omega)\rangle_{L^2(M)} d\mathbb{P}_{\hbar}(\omega)=\frac{1}{N(I(\hbar))}\text{Tr}\left(B\bbbone_{I(\hbar)}(-\hbar^2\Delta)\right).
\end{equation}
\end{itemize}
\end{def1}
We stress the fact that if $(u_{\hbar}(\omega))_{\omega\in(\Omega_{\hbar},\mathbb{P}_{\hbar})}$ is a G.O.F., then $(U^tu_{\hbar}(\omega))_{\omega\in(\Omega_{\hbar},\mathbb{P}_{\hbar})}$ is also one for every $t$. In section~\ref{s:examples}, we will provide two examples of G.O.F.

We will denote $\mu_{\hbar,\omega}(t)$ the (time-dependent) distribution associated to $u_{\hbar}(\omega)$ by formula~(\ref{general-def-distrib}).

\subsubsection{Semiclassical large deviations} We now state a generalization of the quantum ergodicity theorem. This theorem says that, for a given orthonormal basis of eigenvectors of $\Delta$, ``{most of}'' the associated distributions on $T^*M$ converge to the Liouville measure on $S^*M:=\{p=1/2\}$. This holds under the only assumption that the geodesic flow acts ergodically on $S^*M$ endowed with the Liouville measure. Here we aim for a more precise statement, and will assume that the geodesic flow has the Anosov property (our result will, in particular, imply a reinforced version of the usual Shnirelman theorem). 

We recall that the Liouville measure on $T^*M$ is the measure given by $d\cL=dxd\xi$ in local coordinates. In a region where the hamiltonian $p$ has no critical point, one can find local symplectic coordinates $(x_1,\ldots, x_d, \xi_1, \ldots, \xi_d)$ such that $x_1=p$, and the Liouville measure can be decomposed into
$d\cL=dx_1 dL_{x_1}(x, \xi)$, where $L_{x_1} $ is a smooth positive measure carried by the energy layer $\{p=x_1\}$. We shall restrict our attention to the unit cotangent bundle, $S^*M=\{p=\frac12\}$, and
will denote $L=L_{\frac12}$. This is the Liouville measure in $S^*M$.

Given a G.O.F. $(u_\hbar(\omega))_{\omega\in(\Omega_{\hbar},\mathbb{P}_{\hbar})}$, our result says that for ``most'' $\omega$ (in the sense of $\mathbb{P}_{\hbar}$) the distributions $\mu_{\hbar,\omega}(t)$ are close to the Liouville measure $L$. We will use a large deviations result due to Kifer~\cite{Ki} to give an estimate on the proportion of $\omega$ for which $\mu_{\hbar,\omega}(t)$ is far away from $L$. To state our result, we need to introduce two dynamical quantities. First, we define the maximal expansion rate of the geodesic flow on $S^*M$ as
$$\chi_{\max}:=\lim_{t\rightarrow\pm\infty}\frac{1}{t}\log\sup_{\rho\in S^*M}\|d_{\rho}g^t\|.$$
This quantity gives an upper bound on the Lyapunov exponents over $S^*M$ and it is linked to the range of validity of the semiclassical approximation in the Egorov theorem~\cite{BR}. We also introduce, for every $\delta$ in $\mathbb{R}$ and every $a$ in $\mathcal{C}^{\infty}_o(T^*M,\mathbb{R})$ such that $L(a)=0$,
$$H(\delta):=\inf_{s\in \mathbb{R}}\left\{-s\delta+P\left(sa +\varphi^u\right)\right\},$$
where $f\mapsto P(f)$ is the topological pressure of the continuous map $f$ and $\varphi^u$ is the infinitesimal unstable jacobian (see section~\ref{s:dyn-syst} for details). The map $\delta\mapsto -H(\delta)$ is the Legendre transform of $s\mapsto P(sa+\varphi^u)$ which is a smooth and convex function on $\mathbb{R}$. In particular, $-H$ is a convex map on $\mathbb{R}$ and it satisfies $H(0)=0$ and $H(\delta)<0$ for all $\delta\neq 0$ (see \ref{topological-pressure-subsection}).


\begin{theo}\label{deviations-fct-propres-laplacien} Suppose $(S^*M,(g^t))$ satisfies the Anosov property. We fix a generalized orthonormal families $(u_{\hbar}(\omega))_{\omega\in(\Omega_{\hbar},\mathbb{P}_{\hbar})}$ (with $\hbar\rightarrow0^+$). We fix two observables,
\begin{itemize}
\item an element $\theta$ in $L^1(\R,\mathbb{R}_+)$ such that $\int\theta(t)dt=1$,
\item an element $a$ in $\mathcal{C}^{\infty}_o(T^*M,\mathbb{R})$ such that $\int_{S^*M}adL=0$.
\end{itemize}
Then, we have, for any $\delta>0$,
$$\limsup_{\hbar\rightarrow 0}\frac{\log\mathbb{P}_{\hbar}\left(\left\{\omega\in\Omega_{\hbar}: \mu_{\hbar,\omega}( a\otimes\theta)\geq\delta\right\}\right)}{|\log\hbar|}\leq \frac{H(\delta)}{\chi_{\max}}.$$
\end{theo}

From this theorem and the properties of $H(\delta)$, one can deduce the following corollary~:
\begin{coro}\label{coro-deviations-fct-propres-laplacien} Suppose $(S^*M,(g^t))$ satisfies the Anosov property. We fix a G.O.F. $(u_{\hbar}(\omega))_{\omega\in(\Omega_{\hbar},\mathbb{P}_{\hbar})}$ (with $\hbar\rightarrow0^+$). Then, for every $\delta>0$, for every $a\in\mathcal{C}^{\infty}_o(T^*M,\mathbb{C})$, and for every function $\theta$ in $L^1(\R,\mathbb{R}_+)$, we have
\begin{equation}\label{convergence-proba-Shnirelman-poly}
\mathbb{P}_{\hbar}\left(\left\{\omega\in\Omega_{\hbar}: \left|\mu_{\hbar,\omega}(a\otimes\theta)-\int_{S^*M}adL\int_{\mathbb{R}}\theta(t)dt\right|\geq\delta\right\}\right)=\mathcal{O}_{a,\delta,\theta}\left(\hbar^{\tilde{H}(\delta)}\right),
\end{equation}
where $\tilde{H}(\delta)>0$ depends on $a$, $\theta$ and $\delta$.
\end{coro}

\subsubsection{Comments}\label{s:comments-deviations} 
As already mentioned, this result reinforces the Shnirelman theorem in the case of Anosov geodesic flows. The Shnirelman theorem (suitably adapted to the time dependent Schr\"odinger equation) simply asserts that for an ergodic geodesic flow, and for every $\delta>0$,
$$\mathbb{P}_{\hbar}\left(\left\{\omega\in\Omega_{\hbar}: \left|\mu_{\hbar,\omega}(a\otimes\theta)-\int_{S^*M}adL\int_{\mathbb{R}}\theta(t)dt\right|\geq\delta\right\}\right)=o_{a,\delta,\theta}\left(1\right).$$
Let us also recall the quantum variance conjecture~\cite{EFKal}. This conjecture is usually formulated for eigenfunctions of the laplacian, but translated in our context, it would predict that
$$\int_{\Omega_{\hbar}}\left|\mu_{\hbar,\omega}\left(a\otimes\theta\right)-\int_{S^*M}adL\int_{\mathbb{R}}\theta(t)dt\right|^2d\mathbb{P}_{\hbar}(\omega)=\frac{V(a,\theta)\hbar^d}{b(\hbar)-a(\hbar)}(1+o(1))$$
for some constant $V(a,\theta)$ that can be expressed in terms of the ``dynamical variance'' of $a$~\cite{S0}.
If this conjecture is true, it implies
$$\mathbb{P}_{\hbar}\left(\left\{\omega\in\Omega_{\hbar}: \left|\mu_{\hbar,\omega}(a\otimes\theta)-\int_{S^*M}adL\int_{\mathbb{R}}\theta(t)dt\right|\geq\delta\right\}\right)=\mathcal{O}_{a,\delta,\theta}\left(\frac{\hbar^d}{b(\hbar)-a(\hbar)}\right),$$
which is stronger than our result.

Zelditch proved in \cite{Ze2} that
$$\int_{\Omega_{\hbar}}\left|\mu_{\hbar,\omega}\left(a\otimes\theta\right)-\int_{S^*M}adL\int_{\mathbb{R}}\theta(t)dt\right|^pd\mathbb{P}_{\hbar}(\omega)=\cO(|\log\hbar|^{-p/2})$$
for all $p\geq 1$ (see also \cite{S0}). Again, his proof is written for the eigenfunction problem, but could easily be transposed to the time-dependent Schr\"odinger equation (see \cite{Riv1} -- and note that we have to make the extra assumption $\|u_{\hbar}(\omega)\|_{L^2}=1+\mathcal{O}(|\log\hbar|^{-1})$ uniformly in $\omega$).
 Using the Bienaym\'e-Chebyshev inequality, Zelditch's result implies that
$$\mathbb{P}_{\hbar}\left(\left\{\omega\in\Omega_{\hbar}: \left|\mu_{\hbar,\omega}(a\otimes\theta)-\int_{S^*M}adL\int_{\mathbb{R}}\theta(t)dt\right|\geq\delta\right\}\right)=\mathcal{O}(|\log\hbar|^{-\infty}).$$
Our theorem -- although it does not say anything about the quantum variance -- improves this aspect of Zelditch's result, as we can replace $\mathcal{O}(|\log\hbar|^{-\infty})$ by
$\mathcal{O}\left(\hbar^{\tilde{H}(\delta)}\right)$.  

\subsection{Entropy of semiclassical measures}

Our second result is a lower bound on the Kolmogorov-Sinai entropy of the measures $\mu(t)$. We will consider a sequence of normalized states $(u_{\hbar})_{\hbar\rightarrow0^+}$ in $L^2(M)$.
%
We fix two energy levels $0\leq E_1<E_2$ and we suppose that the family of states is localized in the energy window $[E_1,E_2]$. Precisely, we make the assumption that
\begin{equation}\label{e:local-etats}
\lim_{\hbar\rightarrow0^+}\left\|\left(\text{Id}_{L^2(M)}-\bbbone_{[E_1,E_2]}\left(-\hbar^2\Delta\right)\right)u_{\hbar}\right\|_{L^2(M)}=0.
\end{equation}
This assumption implies that each $\mu(t)$ is a probability measure carried by the set $\{E_1\leq\|\xi\|_x^2\leq E_2\}$ (it prevents escape of mass in the fibers of $T^*M$). In addition, we recall that $
\mu(t)$ is invariant under the geodesic flow. Using the invariance of the energy under the geodesic fow, we see that for Lebesgue a.e. $t$, $\mu(t)(dx, d\xi)$ is of the form $\int \mu_{t,E}(dx, d\xi)\nu(dE)$, where $\nu$ is a positive measure on the interval $[E_1,E_2]$ and $\mu_{t,E}$ is a probability measure on $\{\|\xi\|_x^2= E\}$ invariant under the geodesic flow.
\begin{rema} We underline the fact that the measure $\nu$ is independent of $t$. It is the weak limit
(after extraction of a subsequence) of the measures $\nu_\hbar$ defined on $\R$ by $\nu_\hbar([E, E'])=\left\|\bbbone_{[E,E']}\left(-\hbar^2\Delta\right)u_{\hbar}\right\|^2$.
\end{rema}
In the following theorem, $h_{KS}(\mu, (g^t))$ denotes the entropy of the invariant probability measure $\mu$ for the geodesic flow $g^t$ (its definition is recalled in section \ref{s:dyn-syst}).
\begin{theo}\label{t:entropy} Let $M$ be a compact riemannian manifold of dimension $d$ and constant
curvature $\equiv -1$. We fix two energy levels $0\leq E_1<E_2$ and we consider a sequence $(u_{\hbar})_{\hbar\rightarrow0^+}$ in $L^2(M)$ that satisfies:
\begin{itemize}
\item the energy localization
$\lim_{\hbar\rightarrow0}\left\|\left(\text{Id}_{L^2(M)}-\bbbone_{[E_1,E_2]}(-\hbar^2\Delta)\right)u_{\hbar}\right\|_{L^2(M)}=0,$
\item $\lim_{\hbar\rightarrow0}\|u_{\hbar}\|_{L^2(M)}=1.$
\end{itemize}
Consider $\mu(t)(dx, d\xi)=\int \mu_{t,E}(dx, d\xi)\nu(dE) dt$ an accumulation point of the corresponding sequence of distributions $\mu_{\hbar}$ defined by~(\ref{general-def-distrib}). Then, one has, $\text{Leb}\otimes\nu$ almost everywhere,
$$  h_{KS}(\mu_{t,E}, g^t)\geq\frac{d-1}{2}\sqrt{E}.$$
\end{theo}
For the sake of simplicity, we only state and prove the results in the case of constant curvature. The methods from~\cite{A, AN2} for general Anosov manifolds or from~\cite{Riv} for Anosov surfaces could also be adapted in this setting.

\begin{rema}We note that $\sqrt{E}$ is the speed of trajectories of $g^t$ on the energy layer $\{p=\frac{E}{2} \}$. It is also natural to consider the geodesic flow $\phi^t=g^{t/\sqrt{E}}$ parametrized to have speed $1$ on any energy layer, and our results then reads $  h_{KS}(\mu_{t,E}, \phi^{t})\geq\frac{d-1}{2}.$

If one wants, one can avoid assumption \eqref{e:local-etats} and deal with the issue of escape of mass in a different manner~: consider the space $\cS_0$ of smooth functions $a$ on $T^*M$ that are $0$-homogeneous outside a compact set. The distributions $\mu_\hbar$ are bounded in $L^\infty(\R, \cS_0')$, and one can consider convergent subsequences in the corresponding weak-$\ast$ topology. The corresponding limits $\mu\in L^\infty(\R, \cS_0')$ are actually positive for almost all $t$, and each $\mu(t)$ defines a probability measure on $\widehat{T^*M}$, the cotangent bundle compactified by spheres at infinity.  We note that the flow $\phi^t$ can be extended to the spheres at infinity. We can then write $\mu(t)=\int\mu_{t,E}(dx, d\xi)\nu(dE)$ where now $\nu$ is a probability measure on $[0, +\infty]$. Our result reads~:
$  h_{KS}(\mu_{t,E}, g^t)\geq\frac{d-1}{2}\sqrt{E}$ for $0\leq E<+\infty$, and $  h_{KS}(\mu_{t,E}, \phi^t)\geq\frac{d-1}{2}$ for $0< E\leq +\infty$.
\end{rema}

\subsection{Application to controllability}
Once Theorem \ref{t:entropy} is known, it implies the following observability inequality~:
\begin{theo}\label{t:observability} Let $M$ be a compact riemannian manifold of dimension $d$ and constant
curvature $\equiv -1$. Let $a$ be a smooth function on $M$, and define the closed $g^t$-invariant subset of $S^*M$,
$$K_a=\{\rho\in S^*M, a^2(g^t(\rho))=0\,\, \forall t\in\R\}.$$
Assume that the Hausdorff dimension of $K_a$ is $<d$. Then, for all $T>0$, there exists $C_{T, a}>0$ such that, for all $u$~:
\begin{equation}\label{e:observability}\|u\|^2_{L^2(M)}\leq C_{T, a}\int_0^T\| ae^{\imath t\frac\Delta2}u\|^2_{L^2(M)} dt.\end{equation}
\end{theo}
We recall briefly in appendix \ref{s:control} how to deduce Theorem \ref{t:observability} from Theorem \ref{t:entropy}. This follows a classical argument due to Lebeau~\cite{Leb}, who used it to prove the following fact~: if $M$ is an arbitrary riemannian manifold, and if $K_a= \emptyset$, then \eqref{e:observability} holds.

We can give an example where our assumption on $K_a$ holds. Consider a closed geodesic $\gamma$ with a small tubular neigborhood of this geodesic that does not contain another complete geodesic. We take $a$ to be nonzero on the complementary of this neighborhood and $0$ near the closed geodesic. In this case, one has $K_a=\gamma$ so that our condition holds. Another example, in dimension $d=2$, goes as follows~: take a decomposition of the hyperbolic surface
$M$ into ``hyperbolic pairs of pants'' (there are $2g-2$ pairs of pants if $M$ has genus $g$). The boundary of each pair of pants consists of 3 simple closed geodesics. Take a function $a$ supported in a neighbourhood of the union of these $3g-3$ simple closed geodesics, and assume that $a$ does not vanish on the union of these curves. Thus, any geodesic that avoids the support of $a$ must stay inside one of the pairs of pants. If the length of each of the $3g-3$ boundary components is large enough, this will imply that $K_a$ has dimension $<d$, and our condition will be satisfied.
The existence of a hyperbolic pants decomposition with boundary components of
arbitrary large lengths follows, for instance, from Proposition 2.2 in \cite{Rees}. It would be interesting to find a larger variety of geometric situations in which our assumption on $K_a$ holds.

Following the Hilbert uniqueness method, one knowns that inequality \eqref{e:observability} implies the following~: for any $u_0, u_T\in L^2(M)$, for any $T>0$, there exists $f(t, x)\in L^2([0, T]\times M)$ such that the solutions of
$$\imath\frac{\partial u}{\partial t}+\frac\Delta2 u=a(x)f(t, x)$$
with initial condition $u_{|t=0}=u_0$ satisfies $u_{|t=T}=u_T$. This is called the controllability problem.

\subsection*{Organization of the paper}
In section~\ref{s:dyn-syst}, we describe some background in dynamical systems that we will need at different points of the article. In section~\ref{s:examples}, we give two examples of G.O.F. and apply Theorem~\ref{deviations-fct-propres-laplacien} to them. In sections~\ref{s:proof1} and~\ref{s:proof2}, we prove Theorems~\ref{deviations-fct-propres-laplacien} and~\ref{t:entropy}. Finally, in the appendices, we show how to derive an observability result from Theorem~\ref{t:entropy} (appendix~\ref{s:control}) and we give a brief reminder on semiclassical calculus on a manifold (appendix~\ref{pdo-manifold}).

\section{Dynamical systems background}\label{s:dyn-syst}

\subsection{Anosov property}

In this paper, we suppose that $M$ is a smooth, compact, riemannian manifold $M$ of dimension $d$ (without boundary). The geodesic flow on $T^*M$ is the hamiltonian flow associated to the hamiltonian $p(x,\xi)=\frac{\|\xi\|_x^2}2$. We also assume that, for any $E_0>0$, the geodesic flow $g^t$ is Anosov on the energy layer $p^{-1}(\{\frac{E_0}{2}\})\subset T^*M$~: for all $\rho\in p^{-1}(\{\frac{E_0}{2}\})$, we have a decomposition
$$T_{\rho}p^{-1}\left(\left\{E_0/2\right\}\right)=E^u(\rho)\oplus E^s(\rho)\oplus \mathbb{R}X_{p}(\rho),$$
where $X_p$ is the hamiltonian vector field associated to $p$, $E^u$ the unstable space and $E^s$ the stable space~\cite{KaHa}. We can introduce the infinitesimal unstable Jacobian as follows~\cite{BoRue}:
$$\varphi^u(\rho):=-\frac{d}{dt}\left(\det\left(d_{\rho}g^t_{|E^u(\rho)}\right)\right)_{t=0}.$$

\subsection{Kolmogorov-Sinai entropy}
\label{KSentropy}
Let us recall a few facts about Kolmogorov-Sinai (or metric) entropy that can be found for example in~\cite{Wa}. Let $(X,\mathcal{B},T, \mu)$ be a measurable dynamical system, and $\cP:=(P_{\alpha})_{\alpha\in I}$ a finite measurable partition of $X$, i.e. a finite collection of measurable subsets that forms a partition. Each $P_{\alpha}$ is called an atom of the partition. With the convention $0\log 0=0$, one defines
\begin{equation}\label{nKSentropy}H_{n}(\mu, T, \cP)=-\sum_{|\alpha|=n}\mu(P_{\alpha_0}\cap\cdots\cap T^{-(n-1)}P_{\alpha_{n-1}})\log\mu(P_{\alpha_0}\cap\cdots\cap T^{-(n-1)}P_{\alpha_{n-1}}).\end{equation}
This quantity satisfies a subadditivity property
\begin{equation}\label{e:subadd-KS-entropy}H_{n+m}(\mu,T,\mathcal{P})\leq H_n(\mu,T,\mathcal{P})+H_m\left(\mu,T,T^{-n}\mathcal{P}\right)=H_n(\mu,T,\mathcal{P})+H_m(\mu,T,\mathcal{P}).\end{equation}
The first inequality is true even if the probability measure $\mu$ is not $T$-invariant, while the second equality holds for $T$-invariant measures. A classical argument for subadditive sequences allows to define the following quantity:
\begin{equation}\label{defentropy}h_{KS}(\mu,T,\cP):=\lim_{n\rightarrow\infty}\frac{H_n\left(\mu,T,\cP\right)}{n}.\end{equation}
It is called the Kolmogorov Sinai entropy of $(T,\mu)$ with respect to the partition $P$. The Kolmogorov Sinai entropy $h_{KS}(\mu,T)$ of $(\mu,T)$ is then defined as the supremum of $h_{KS}(\mu,T,\cP)$ over all partitions $\cP$ of $X$. In the case of a flow (for instance the dynamical system $(S^*M,g^t,\mu)$), we define the entropy $h_{KS}(\mu,(g^t)):=h_{KS}(\mu,g^1)$. Entropy can {\em a priori} be infinite. However, for a smooth flow on a compact finite dimensional manifold, entropy is bounded thanks to the Ruelle inequality~\cite{R}
$$h_{KS}(\mu,(g^t))\leq -\int_{S^*M}\varphi^u(\rho)d\mu(\rho).$$
In the case of negatively curved manifolds, equality holds if and only if $\mu$ is the desintegration $L$ of the Liouville measure on $S^*M$. 

{\bf Notation~:} In the rest of this paper, we will write $h_{KS}(\mu)$ for $h_{KS}(\mu,(g^t))$, unless we want to consider a flow different from $(g^t)$.


\subsection{Topological pressure}\label{topological-pressure-subsection} To conclude this section, we introduce the topological pressure of the dynamical system $(S^*M,g^t)$ as the Legendre transform of the Kolmogorov-Sinai entropy~\cite{Wa},~\cite{PaPo}:
$$\forall f\in\mathcal{C}^0(S^*M,\mathbb{R}),\ P(f)=P(f, (g^t)):=\sup\left\{h_{KS}(\mu)+\int_{S^*M}fd\mu:\mu\in\mathcal{M}(S^*M,g^t)\right\},$$
where $\mathcal{M}(S^*M,g^t)$ is the set of probability measures on $S^*M$ invariant under the geodesic flow. It can be verified that this defines a continuous and convex function on $\mathcal{C}^0(S^*M,\mathbb{R})$~\cite{Wa}.

We shall be particularly interested in the behaviour of $P(f)$ near $f=\varphi^u$. By the Ruelle inequality, we have $P(\varphi^u)=0$ (the $\sup$ defining $P(\varphi^u)$ is achieved at $\mu=L$).
Moreover, it can be proved that for any real-valued H\"older function $f$ on $S^*M$, the function $s\mapsto P(\varphi^u+sf)$ is real analytic on $\R$~\cite{BoRue, R0} and its derivatives of order $1$ and $2$ can be computed explicitly~\cite{PaPo}. We have $\frac{d}{d s}\left(P( \varphi^u+sf)\right)_{|s=0}=\int_{S^*M}f dL$.
If $\int_{S^*M}f dL=0$, the convex function $s\mapsto P(\varphi^u+sf)$ achieves its minimum at $0$.

Moreover, if $\int_{S^*M}fdL=0$, then we have:
$$\frac{d^2}{d s^2}\left(P(  \varphi^u+sf)\right)_{|s=0}=\sigma^2(f),$$
where $\displaystyle\sigma^2(f):=\lim_{T\rightarrow+\infty}\frac{1}{T}\int_{S^*M}\left(\int_0^Tf\circ g^t(\rho)dt\right)^2dL(\rho)$ is called the dynamical variance of the function $f$. It is known that $\sigma^2(f)$ vanishes if and only if $f$ is of the form $f=\frac{d}{dt}(h\circ g^t)_{|t=0}$ for some function $h$. In this case, one says that $f$ is a coboundary.

\subsection{Kifer's large deviation upper bound}
We shall use the following result, due to Kifer \cite{Ki}, and valid for more general Anosov flows~:
\begin{equation}\lim_{T\rightarrow\infty}\frac{1}{T}\log\int_{S^*M}\exp\left(\int_{0}^{T}a\circ g^t(\rho)dt\right)dL(\rho)=P(a+\varphi^u),\label{e:kif}
\end{equation}
for all continuous $a$. In fact, we will only use the fact that the $\limsup$ is uniform for $a$ running
over compact sets (section~$3$ -Theorem $3.2$ in~\cite{Ki}).

\begin{rema} This result can be used to proved the following strengthened version of the Birkhoff ergodic theorem. Fix $a$ such that $\int_{S^*M}adL=0$, and fix $\delta>0$. Then,
\begin{multline*}\limsup\frac1T\log L\left(\left\{\rho\in S^*M: \frac1T\int_{0}^{T}a\circ g^t(\rho)dt >\delta\right\}\right)\leq
\inf_{s\geq 0}\{-s\delta +P(sa+\varphi^u)\}\\
=\inf_{s\in \R}\{-s\delta +P(sa+\varphi^u)\}=H(\delta).
\end{multline*}
Similarly, for $\delta<0$, one has $\limsup\frac1T\log L(\{\rho\in S^*M: \frac1T\int_{0}^{T}a\circ g^t(\rho)dt <\delta\})\leq H(\delta)$.
The function $-H$, which is the Legendre transform of $s\mapsto P(\varphi^u+sa)$, is convex and is positive for $\delta\neq 0$ (it may be infinite).
\end{rema}

\section{Examples of generalized orthonormal families}\label{s:examples}

In this section, we provide two examples of G.O.F. and show how Theorem~\ref{deviations-fct-propres-laplacien} applies to them. Our examples are of distinct types: basis of eigenvectors of $\Delta$ and truncated Dirac distributions. In the first example, Theorem~\ref{deviations-fct-propres-laplacien} provides a strengthened version of Shnirelman's theorem for Anosov flows.

\subsection{Orthonormal basis of eigenvectors}

Consider $(\psi_n)_{n\in\mathbb{N}}$ an orthonormal basis of $L^2(M)$ made of eigenfunctions of $\Delta$, i.e. there exists a sequence $0=\lambda_0<\lambda_1\leq\cdots\leq\lambda_n\leq\cdots$ such that for every $n$ in $\mathbb{N}$,
$$\Delta\psi_n=-\lambda_n^2\psi_n.$$
For $\hbar>0$,we take $\Omega_{\hbar}:=\{n\in\mathbb{N}:\hbar^2\lambda_n^2\in[1-\alpha\hbar,1+\alpha\hbar]\},$ where $\alpha$ is some fixed positive number. In this case, the probability measure is given by $\mathbb{P}_{\hbar}:=\frac{1}{|\Omega_{\hbar}|}\sum_{n\in\Omega_{\hbar}}\delta_n$ and the measurable map is given by $u_{\hbar}(n):=\psi_n$. Applying corollary~\ref{coro-deviations-fct-propres-laplacien} to this example, we find that for every $a$ in $\mathcal{C}^{\infty}_o(T^*M)$, and for every $\delta>0$, there exists $\tilde{H}(\delta)>0$ such that
$$\frac{1}{|\Omega_{\hbar}|}\left|\left\{n\in\Omega_{\hbar}:\left|\mu_{\hbar,n}(a)-\int_{S^*M}adL\right|\geq\delta\right\}\right|=\mathcal{O}_{a,\delta}(\hbar^{\tilde{H}(\delta)}).$$
Shnirelman theorem provides a $o_{a, \delta}(1)$ and using the results from~\cite{Ze2} on eigenfunctions of $\Delta$, one would obtain a $\mathcal{O}_{a, \delta, p}(|\log\hbar|^{-p})$ for arbitrarily large $p$.

\subsection{Truncated Dirac distributions}
The second class of examples we will consider is given by families of vectors constructed from the Dirac distributions. For $y$ in $M$, we denote $\delta_y$ the Dirac distribution given by $\langle\delta_y,f\rangle:=f(y)$ (where $f$ is in $\mathcal{C}^{\infty}(M)$). To construct our G.O.F., we will project $\delta_y$ on $L^2(M)$. To do this, recall that we have defined $I(\hbar):=[a(\hbar),b(\hbar)]$, where $b(\hbar)-a(\hbar)\geq 2\alpha\hbar$ and that we have denoted $N(I(\hbar)):=|\{n:\hbar^2\lambda_n^2\in I(\hbar)\}|$. Using these notations, we can introduce a truncated Dirac distribution as follows:
$$\delta_y^{\hbar}:=\left(\frac{\text{Vol}_M(M)}{N(I(\hbar))}\right)^{\frac{1}{2}}\bbbone_{I(\hbar)}\left(-\hbar^2\Delta\right)\delta_y.$$
According to (global and local) Weyl laws from~\cite{DuisGui} and from~\cite{SogZe} (Theorem $1.2$), we know that in the Anosov case,
$$\left(M,\frac{\text{Vol}_M}{\text{Vol}_M(M)},\delta_y^{\hbar}\right)\ \text{is a G.O.F. in the spectral window}\ I(\hbar).$$
Applying Corollary~\ref{coro-deviations-fct-propres-laplacien} to this example, we find that for every $a$ in $\mathcal{C}^{\infty}_o(T^*M)$, for every $\theta$ in $L^1(\mathbb{R},\mathbb{R}_+)$ and for every $\delta>0$, there exists $\tilde{H}(\delta)>0$ such that
$$\text{Vol}_M\left(\left\{y\in M:\left|\mu_{\hbar,y}(a\otimes\theta)-\int_{S^*M}adL\int_{\mathbb{R}}\theta(t)dt\right|\geq\delta\right\}\right):=\mathcal{O}_{a,\theta,\delta}(\hbar^{\tilde{H}(\delta)}).$$
Thus, if we choose $y$ randomly on $M$ according to the volume measure, and consider the solution
of the Schr\"odinger equation $e^{\imath t\frac{\Delta}2}\delta_y^{\hbar}$,
our result says that we have convergence of the associated semiclassical measure to the uniform  measure, for most $y$ (in the probability sense, and with an explicit bound) as $\hbar$ tends to $0$. Taking a subsequence $(\hbar_n)_n$ that tends to $0$ fast enough, we can apply the Borel-Cantelli lemma and derive convergence for almost every $y$~\cite{Riv1}. An interesting question would be to understand more precisely for which subsequences $(\hbar_n)$ we have convergence for almost every $y$.

\subsection{Coherent states} Similar results could, in principle, apply to bases of coherent states (e.g. gaussian states). Such bases can be constructed easily in euclidean situations; see \cite{Riv1} for an application of Theorem \ref{deviations-fct-propres-laplacien} to the ``cat-map'' toy model. However, on an arbitrary manifold, it seems difficult to construct bases of coherent states meeting all the requirements of the definition of a G.O.F, which are actually quite strong.

\section{Proof of theorem~\ref{deviations-fct-propres-laplacien}}\label{s:proof1}

In this section, we give a proof of Theorem~\ref{deviations-fct-propres-laplacien}: there are two steps. To begin with, we combine the Bienaym\'e-Chebyshev inequality and the Egorov theorem to obtain a first bound (\S \ref{s:BT}). Then we apply a large deviations estimate due to Kifer~\cite{Ki} to obtain a bound in terms of the topological pressure.\\
We fix $\theta$ an element of $L^1(\R,\mathbb{R}_+)$ such that $\int \theta(t)dt=1$. Let $a$ be an element in $\mathcal{C}^{\infty}_o(T^*M,\mathbb{R})$ that satisfies $\int_{S^*M}adL=0$. Recall that we defined
$$\chi_{\max}:=\lim_{t\rightarrow\pm\infty}\frac{1}{t}\log\sup_{\rho\in S^*M}\|d_{\rho}g^t\|.$$
As the states $u_{\hbar}(\omega)$ are uniformly microlocalized in a thin neighborhood of $S^*M$, we can assume that $a$ is compactly supported in a tubular neighborhood $p^{-1}\left([\frac12-\eta, \frac12+\eta]\right)$ of $S^*M$ (with $\eta>0$ arbitrarily small). Letting $\chi_\eta=\chi_{\max}\sqrt{1+2\eta}$, we have
$$\forall t\in\mathbb{R},\ \forall \rho\in T^{*}M, \forall \alpha, \ \|\partial^{\alpha}(a\circ g^t)(\rho)\|\leq C_{a,\alpha} e^{\chi_\eta|\alpha| |t|}.$$

\subsection{Long-time Egorov theorem}
We fix $c$ such that $c\chi_\eta<\frac12$. 
The positive quantization $\Op_{\hbar}^+$ procedure satisfies the following ``long time Egorov property''~:
\begin{equation}\forall |t|\leq c|\log\hbar|,\ \|U^{-t\hbar}\Op_{\hbar}^+(a)U^{t\hbar}-\Op_{\hbar}^+(a\circ g^t)\|_{L^2(M)\rightarrow L^2(M)}=\mathcal{O}_a(\hbar^{\frac{1}{2}-\nu}),\label{e:egorov}
\end{equation}
where $\nu:=c\chi_\eta$ (see \cite{AN2}).
\begin{lemm}\label{invariance-lemma} Using previous notations, for every $\delta_0>0$, there exists $\hbar_0$ (depending on $a$, $\theta$ and $\delta_0$) such that for every $\hbar<\hbar_0$, we have for every $|T|\leq c|\log\hbar|$:
$$\left\|\int\theta(t)U^{-t}\left(\Op_{\hbar}^+(a)-\frac{1}{2T}\int_{-T}^T\Op_{\hbar}^+(a\circ g^s)ds\right)U^tdt\right\|_{L^2(M)\rightarrow L^2(M)}\leq\delta_0.$$
\end{lemm}
\emph{Proof.} The proof of this lemma relies on the application of the Egorov property \eqref{e:egorov}. For  $T$ a real number such that $|T|\leq c|\log\hbar|$, we have$$\int\theta(t)U^{-t}\left(\frac{1}{2T}\int_{-T}^T\Op_{\hbar}^+(a\circ g^s)ds\right)U^tdt=\frac{1}{2T}\int_{-T}^T\int\theta(t)U^{-t-s\hbar}\Op_{\hbar}^+(a)U^{t+s\hbar}dtds+\mathcal{O}_a(\hbar^{\frac{1}{2}-\nu}).$$
We make the change of variables $t'=t+s\hbar$ and use the fact that $\|\theta(.)-\theta(.-\tau))\|_{L^1}\Lim_{\tau\To 0}0$ to conclude.$\square$

\subsection{Bienaym\'e-Chebyshev and Jensen's inequality\label{s:BT}}

For simplicity of notations, we will denote the quantity we want to bound as follows:
$$\mathbb{P}_{\hbar}(a\otimes\theta,\delta):=\mathbb{P}_{\hbar}\left(\left\{\omega\in\Omega_{\hbar}: \mu_{\hbar,\omega}( a\otimes\theta)\geq\delta\right\}\right).$$

 Our first step is to combine the previous lemma to the Bienaym\'e-Chebyshev inequality in order to obtain a bound on $\mathbb{P}_{\hbar}(a\otimes\theta,\delta)$.
\begin{lemm}\label{etape-1-preuve-deviations-M} Let $\delta, \delta_0>0$ be arbitrary positive numbers. For $s\in\R$, denote
$$a_{s}(T(\hbar),\rho):=\exp\left(s\int_{-T(\hbar)}^{T(\hbar)}a\circ g^\tau(\rho)d\tau\right),$$
where $T(\hbar)=c|\log\hbar|$ (and $c$ is such that $c\chi_\eta<1/2$). Then, given $s>0$ and for $\hbar$ small enough, one has
\begin{equation}  \mathbb{P}_{\hbar}(a\otimes\theta,\delta)\leq 2 \frac{e^{(-2\delta+3\delta_0)s T(\hbar)}}{N(I(\hbar))}\Tr\left[\bbbone_{I(\hbar)}(-\hbar^2\Delta ) \Op_{\hbar}^+\left(a_{s}(T(\hbar),\bullet)\right) \right].\end{equation}
\end{lemm}
\emph{Proof.} To prove this lemma, we fix $s>0$. A direct application of the Bienaym\'e-Chebyshev inequality allows us to write
$$\mathbb{P}_{\hbar}(a\otimes\theta,\delta):=\mathbb{P}_{\hbar}\left(\left\{\omega\in\Omega_{\hbar}: \mu_{\hbar,\omega}(a\otimes\theta )\geq\delta\right\}\right)\leq e^{-2s\delta T(\hbar)}\int_{\Omega_{\hbar}}\exp\left(2sT(\hbar)\mu_{\hbar,\omega}(a\otimes\theta)\right)d\mathbb{P}_{\hbar}(\omega).$$
We can now use Lemma~\ref{invariance-lemma} to deduce that, for $\hbar$ small enough,
$$\mathbb{P}_{\hbar}(a\otimes\theta,\delta)\leq e^{-2s\delta T(\hbar)}\int_{\Omega_{\hbar}}\exp\left(s\mu_{\hbar,\omega}\left(\left(\int_{-T(\hbar)}^{T(\hbar)}a\circ g^\tau d\tau\right)\otimes\theta\right)+2s\delta_0T(\hbar)\|u_{\hbar}(\omega)\|^2\right)d\mathbb{P}_{\hbar}(\omega).$$
Using the fact that $\|u_{\hbar}(\omega)\|=1+o(1)$ uniformly for $\omega$ in $\Omega_{\hbar}$, the quantity $e^{2s\delta_0T(\hbar)\|u_{\hbar}(\omega)\|}$ is uniformly bounded by $e^{3s\delta_0T(\hbar)}$ for $\hbar$ small enough. The map $x\mapsto e^{sx}$ is convex and we can use Jensen's inequality to write
$$\mathbb{P}_{\hbar}(a\otimes\theta,\delta)\leq e^{s(-2\delta+3\delta_0) T(\hbar)}\int_{\Omega_{\hbar}}\mu_{\hbar,\omega}\left(\exp\left(s\mu_{\hbar,\omega}(1\otimes\theta)\left(\int_{-T(\hbar)}^{T(\hbar)}a\circ g^\tau d\tau\right)\right)\otimes\theta\right)\frac{d\mathbb{P}_{\hbar}(\omega)}{\mu_{\hbar,\omega}(1\otimes\theta)}.$$
Using again that $\|u_{\hbar}(\omega)\|=1+o(1)$ uniformly for $\omega$ in $\Omega_{\hbar}$, we find (using that $\theta$ is nonnegative and of integral $1$) $$0<\frac{1}{2}\leq\mu_{\hbar,\omega}(1\otimes\theta)\leq1+\frac{\delta_0}{\|a\|_{\infty}},$$ uniformly for $\hbar$ small enough. All this can be summarized as follows:
$$\mathbb{P}_{\hbar}(a\otimes\theta,\delta)\leq2e^{s(-2\delta+4\delta_0) T(\hbar)}\int_{\Omega_{\hbar}}\mu_{\hbar,\omega}\left(a_{s}(T(\hbar),\bullet)\otimes\theta\right)d\mathbb{P}_{\hbar}(\omega).$$
Note that the function $a_{s}(T(\hbar),\bullet)$ belongs to the class of symbols $S^{0,k_0}_{\nu}(T^*M)$ where $\nu:=c\chi_{\eta}<1/2$ and $k_0:=2cs\|a\|_{\infty}$; moreover $a_s(T(\hbar),\bullet)$ is constant in a neighborhood of infinity. The previous inequality can be rewritten as~:
$$\mathbb{P}_{\hbar}(a\otimes\theta,\delta)\leq 2e^{(-2\delta+4\delta_0)s T(\hbar)}\int\theta(t)\int_{\Omega_{\hbar}}\langle u_{\hbar}(\omega)|U^{-t}\Op_{\hbar}^+(a_s(T(\hbar),\bullet))U^t|u_{\hbar}(\omega)\rangle d\mathbb{P}_{\hbar}(\omega)dt.$$
We recall that if $(u_{\hbar}(\omega))_{\omega\in(\Omega_{\hbar},\mathbb{P}_{\hbar})}$ is a G.O.F. then for every $t$ in $\mathbb{R}$, $(U^tu_{\hbar}(\omega))_{\omega\in(\Omega_{\hbar},\mathbb{P}_{\hbar})}$ is also a G.O.F. Using point $3$ of the definition of a G.O.F, we get the following bound for $\hbar$ small enough:
\begin{equation}\mathbb{P}_{\hbar}(a\otimes\theta,\delta)\leq \frac{2e^{(-2\delta+4\delta_0)s T(\hbar)}}{N(I(\hbar))}\Tr\left[\bbbone_{I(\hbar)}(-\hbar^2\Delta ) \Op_{\hbar}^+\left(a_{s}(T(\hbar),\bullet)\right) \right].\end{equation}
$\square$

\subsection{Trace asymptotics} We now have to estimate (from above) the trace\begin{equation}\label{e:trace}\Tr\left[\bbbone_{I(\hbar)}(-\hbar^2\Delta ) \Op_{\hbar}^+\left(a_{s}(T(\hbar),\bullet)\right) \right] .
\end{equation}
We first underline that, for every $\hbar>0$, there exist energy levels $E_1<\cdots<E_P$ (depending on $\hbar$) such that
$$I(\hbar)=[a(\hbar),b(\hbar)]\subset\bigsqcup_{p=1}^P[E_p-\alpha\hbar,E_p+\alpha\hbar)\subset[a(\hbar)-\alpha\hbar,b(\hbar)+\alpha\hbar],$$
for some fixed positive $\alpha$. Note that $P=\mathcal{O}((b(\hbar)-a(\hbar))/\hbar)$.
We decompose \eqref{e:trace} into
$$\sum_{p=1}^P\Tr\left[\bbbone_{[E_p-\alpha\hbar,E_p+\alpha\hbar)}(-\hbar^2\Delta ) \Op_{\hbar}^+\left(a_{s}(T(\hbar),\bullet)\right) \right] .$$
We shall bound each term of the previous sum (uniformly with respect to $p$), using standard trace estimates, and then sum over $p$. We consider for instance the interval $[1-\alpha\hbar,1+\alpha\hbar)$, and recall how to determine the asymptotic behaviour of 
$$ \Tr\left[\bbbone_{[1-\alpha\hbar,1+\alpha\hbar)}(-\hbar^2\Delta ) \Op_{\hbar}^+\left(a_{s}(T(\hbar),\bullet)\right) \right].$$
Introduce a function $f$ which is $\mathcal{C}^{\infty}$, compactly supported in a small neighborhood of $1$, equal to $1$ in a neighbourhood of $1$ and taking values in $[0,1]$. We shall also use a function $\chi$ in $\mathcal{S}(\mathbb{R}^d)$ whose Fourier transform is compactly supported in a small neighborhood of $0$, containing no period of the closed geodesics of $(g^t)$ on $S^*M$. We assume that $\chi\geq0$ and that it is greater than $1$ on $[-\alpha,\alpha]$. Using the fact that the quantization is positive, we can bound the previous quantity as follows:
\begin{multline}\label{etape-positif-weyl-preuve}\Tr\left[\bbbone_{[1-\alpha\hbar,1+\alpha\hbar)}(-\hbar^2\Delta ) \Op_{\hbar}^+\left(a_{s}(T(\hbar),\bullet)\right) \right]\\
\leq\Tr\left[f(-\hbar^2\Delta )\chi\left(\frac{-\hbar^2\Delta -1}\hbar\right) \Op_{\hbar}^+\left(a_{s}(T(\hbar),\bullet)\right) \right] . 
\end{multline}
The study of this last quantity now follows well known lines. We use the Fourier inversion formula, $\displaystyle2\pi\chi\left(\frac{E-1}{\hbar}\right)=\int_{\mathbb{R}} e^{\frac{\imath(E-1)}{\hbar}t}\hat{\chi}(t)dt$. As a consequence, the right-hand side of \eqref{etape-positif-weyl-preuve} can be written as
$$\frac{1}{2\pi}\int_{\mathbb{R}}e^{-\frac{\imath}{\hbar}t}\text{Tr}\left(\Op_{\hbar}^+(a_{s}(T(\hbar),\bullet))U^{2t\hbar}f(-\hbar^2\Delta)\right) \hat{\chi}(t)dt.$$
The asymptotic behaviour of the trace comes from an asymptotic expansion of the kernel of the operator $\Op_{\hbar}^+(a_{s}(T(\hbar),\bullet))U^{2t\hbar}f(-\hbar^2\Delta)$. This expansion is given by the theory of Fourier integral operators~\cite{DS} (chapter $11$),~\cite{EZ} (chapter $10$). The trace is then expressed as the integral of the kernel over the diagonal, and the asymptotic behaviour of this integral is determined thanks to the method of stationary phase (\cite{DS}, chapter $11$).

 \begin{lemm}\label{prop-semi-classique-schubert} For every integer $N\geq 1$, we have \begin{multline*}\Tr\left[f(-\hbar^2\Delta )\chi\left(\frac{-\hbar^2\Delta -1}\hbar\right) \Op_{\hbar}^+\left(a_{s}(T(\hbar),\bullet)\right) \right]=
 \\  \frac{ 1}{(2\pi\hbar)^{d-1}}\left(\sum_{n=0}^{N-1}\hbar^n\int_{S^*M}D^{2n}a_{s}(T(\hbar),\rho)dL(\rho)+\cO_{a, \chi, \theta, N}(\hbar^{N(1-2\nu)-\beta\nu-k_0})\right),
 \end{multline*}
where $\beta>0$ depends only on the dimension of $M$, and where $D^{2n}$ is a differential operator of order $2n $ on $T^*M$ (depending on the cutoff functions and on the choice of the quantization $\Op_\hbar^+$).
\end{lemm}

There are many references for this kind of estimates. For instance, a very similar calculation is done by Schubert in~\cite{S0} (proposition $1$; he stops at $N=1$ but the stationary phase method actually provides
asymptotic expansions at any order).

Recall that $\nu=c\chi_\eta<\frac12$. It is important here to note that $a_{s}(T(\hbar),\bullet)$ belongs to the class $S^{0,k_0}_{\nu}(T^*M)$. We also underline the fact that the observable $a_{s}(T(\hbar),x,\xi)$ satisfies the particular property that $D^{2n}a_{s}(T(\hbar),\rho)$ is of the form $a_{s}(T(\hbar),x,\xi)b_{2n}(x,\xi)$, with $\|b_{2n}\|_{\infty}=\cO(|s|^{2n}\hbar^{-2n\nu})$. If $s$ stays in a bounded interval, and if we choose $N$ large enough accordingly, this implies that
\begin{multline*}\Tr\left[f(-\hbar^2\Delta )\chi\left(\frac{-\hbar^2\Delta -1}\hbar\right) \Op_{\hbar}^+\left(a_{s}(T(\hbar),\bullet)\right) \right]\leq
 \\  \frac{ 1}{(2\pi\hbar)^{d-1}}\left(\int_{S^*M}a_{s}(T(\hbar),\rho)dL(\rho) \right)(1+\cO(\hbar^{1-2\nu})).
 \end{multline*}

Combing this with Lemma~\ref{etape-1-preuve-deviations-M}, using the Weyl law \eqref{e:weyl}, we finally have, for every $N\geq 1$ and $\hbar$ small enough,
\begin{equation}\label{conclusion-bienayme-M}\mathbb{P}_{\hbar}(a\otimes\theta,\delta)\leq Ce^{(-2\delta+4\delta_0)s T(\hbar)}\left(\int_{S^*M}a_{s}(T(\hbar),\rho)dL(\rho)\right)(1+\cO(\hbar^{1-2\nu})),
\end{equation}
for some constant $C$ that does not depend on $\hbar$.

\subsection{A large deviations bound}

To conclude, we use Kifer's large deviations result \eqref{e:kif}. For our proof, we only need an upper bound on the quantity 
$$\int_{S^*M}\exp\left(s\int_{-T}^{T}a\circ g^t(\rho)\right)dL(\rho).$$
Compared with \eqref{e:kif}, there is a parameter $s$ in the exponential that stays in a bounded interval~$I$. We use the upper bound \eqref{e:kif}, which shows that for every $\delta'>0$ and any bounded interval $I$ in $\mathbb{R}_+$, there exists $c_{\delta'}>0$ and $n(\delta',I)\in\mathbb{N}$ such that for every $n\geq n(\delta',I)$ and every $s$ in $I$:
\begin{equation}\label{bound-largedeviations-flots}\int_{S^*M}\exp\left(s\int_{-T}^{T}a\circ g^t(\rho)\right)dL(\rho)\leq c_{\delta'}e^{T\delta'}e^{2T P\left(s a+\varphi^u\right)}.\end{equation}
This last bound will allow us to conclude. In fact, combining this inequality to the bound~(\ref{conclusion-bienayme-M}) on $\mathbb{P}_{\hbar}(a\otimes\theta,\delta)$, we find that:
$$\mathbb{P}_{\hbar}(a\otimes\theta,\delta)\leq C e^{(-2\delta+4\delta_0)s T(\hbar)} e^{T(\hbar)\delta'}e^{2T(\hbar) P\left(s a+\varphi^u\right)} ,$$
where the constant $C$ depends on the different parameters but not on  $\hbar$. This implies
$$\limsup_{\hbar\rightarrow0}\frac{\log\left(\mathbb{P}_{\hbar}(a\otimes\theta,\delta)\right)}{c|\log\hbar|}\leq \delta'+(-2\delta+4\delta_0)s+2P\left(s a+\varphi^u\right).$$
This last inequality holds for any $\delta_0>0$ and any $\delta'>0$. It implies that for every $s>0$ in the interval $I$:
$$\forall c\in\left(0,\frac{1}{2\chi_{\max}}\right) ,\ \limsup_{\hbar\rightarrow0}\frac{\log\left(\mathbb{P}_{\hbar}(a\otimes\theta,\delta)\right)}{c|\log\hbar|}\leq -2s\delta+2P\left(s a+\varphi^u\right).$$
In particular, we find that
$$\forall\delta\in\mathbb{R},\ \limsup_{\hbar\rightarrow0}\frac{\log\left(\mathbb{P}_{\hbar}(a\otimes\theta,\delta)\right)}{\frac{|\log\hbar|}{2\chi_{\max}}}\leq 2 \inf_{s\in \mathbb{R}_+}\left\{-s\delta+P\left(s a+\varphi^u\right)\right\}.$$
Since $\delta>0$, we have
$\inf_{s\in \mathbb{R}_+}\left\{-s\delta+P\left(s a+\varphi^u\right)\right\}=\inf_{s\in \mathbb{R}}\left\{-s\delta+P\left(s a+\varphi^u\right)\right\}.$
This concludes the proof of Theorem~\ref{deviations-fct-propres-laplacien}.$\square$

\section{Proof of Theorem~\ref{t:entropy}}\label{s:proof2}

In this section, we fix two energy levels $0\leq E_1<E_2$ and consider a sequence $(u_{\hbar})_{\hbar\rightarrow0^+}$ in $L^2(M)$ that satisfies 
$$\lim_{\hbar\rightarrow0}\left\|\left(\text{Id}_{L^2(M)}-\bbbone_{[E_1,E_2]}(-\hbar^2\Delta)\right)u_{\hbar}\right\|_{L^2(M)}=0.$$
Moreover, we suppose that $ \|u_{\hbar}\|_{L^2(M)}=1.$ 
The proof follows essentially the same lines as in~\cite{AN2}, and we refer the reader to that paper for a detailed account. 

\subsection{Quantum partitions}\label{s:smooth-partition}

As usual when computing the Kolmogorov--Sinai entropy, we start by decomposing the manifold $M$
into finitely many pieces (of small diameter). Let $(P_k)_{k=1,\ldots,K}$ be a family of smooth real functions on $M$,
with $\text{supp} P_k\Subset \Omega_k$, such that
\begin{equation}\label{e:partition}
\forall x\in M,\qquad \sum_{k=1}^K P_k^2(x)= 1\,.
\end{equation}
Later on we will assume that the diameters of the supports of the $P_k$ are small enough.
We shall denote $\hat P_k$
the operator of multiplication by $P_k(x)$ on the Hilbert space $L^2(M)$. We denote the Schr\"odinger flow by $U^t=\exp(\frac{\imath t\Delta}{2})$. With no loss of generality, we
will assume that
the injectivity radius of $M$ is greater than $2$, and work with this
propagator at time $\hbar$, $U^{\hbar}$. This unitary operator is a Fourier integral operator associated with the geodesic flow taken at time $t=1$, $g^1$. As one does to compute the Kolmogorov-Sinai entropy of an invariant measure,
we define a new quantum partition of unity by evolving and refining
the initial partition under
the quantum evolution. For each time $n\in\mathbb{N}$ and any
sequence of symbols $\alpha=(\alpha_0,\cdots,\alpha_{n-1})$,
$\alpha_i\in [1,K]$ (we say that the sequence ${\alpha}$ is of {\em length} $|{\alpha}|=n$),
we define the operators
\begin{equation}\label{e:P_bep}
\pi_{\alpha}= \hat{P}_{\alpha_{n-1}}(n-1)\hat{P}_{\alpha_{n-2}}(n-2)\ldots \hat{P}_{\alpha_0}\,.
\end{equation}
Throughout the paper we use the notation $\hat{A}(t)=U^{-t\hbar} \hat A U^{t\hbar}$ for the quantum evolution of an operator $\hat A$. From \eqref{e:partition} and the unitarity of $U$,
the family of operators $\{\pi_{\alpha}:|\alpha|=n\}$ obviously satisfies the resolution of identity
$\sum_{|{\alpha}|=n} \pi_{\alpha}\pi_{\alpha}^{*}=\text{Id}_{L^2(M)}$. We also have $\sum_{|\alpha|=n} \pi^*_{\alpha} \pi_{\alpha}=\text{Id}_{L^2(M)}$.

\subsection{Quantum entropy, and entropic uncertainty principle}\label{s:mainapp}
 For each time $n$, and each normalized $\phi\in L^2(M)$, we define two quantities that 
 are noncommutative analogues of the entropy \eqref{nKSentropy}~:
 \begin{align}\label{e:entropy}
h^-_n(\phi)&=-\sum_{|\alpha|=n}\left\|\pi^*_{\alpha}\phi\right\|^2\log\left(\left\|\pi^*_{\alpha}\phi\right\|^2\right),\\
h^+_n(\phi)&=-\sum_{|\alpha|=n}\left\|\pi_{\alpha}\phi\right\|^2\log\left(\left\|\pi_{\alpha}\phi\right\|^2\right)
\end{align}
In all that follows, the integer $n$ is of order $\tilde c|\log\hbar|$ (with $\tilde c>0$ to be chosen later), and thus the
number of terms in the sum $\sum_{|\alpha|=n}$ is of order $\hbar^{- K_0}$ for some $K_0>0$.
The following is proved in \cite{AN2}, using the entropic uncertainty principle of \cite{MU}.
\begin{prop}\label{c:WEUP}
Let $\chi$ be real-valued, smooth, compactly supported function on $\R$.
Define
\begin{equation}\label{e:c_C}
c(\chi, n):=
\max_{|\alpha|=|\alpha'|=n}\left(
\|\pi_{\alpha'}(n)\pi_{\alpha}\, \chi(-\hbar^2\Delta)\|\right).
\end{equation}
Then for any $\hbar>0$, $L>0$, for any normalized state $\phi$ satisfying
,
\begin{equation}\label{e:local}\sup_{|\alpha|=n}\|(I- \chi(-\hbar^2\Delta))\pi^*_\alpha\phi\|\leq \hbar^L,
\end{equation}
we have
$$h_n^+( U^{n\hbar}\phi) + h_n^-( \phi) \geq
-2\log\left(c(\chi, n)+h^{L-K_0}\right).
$$
\end{prop}

Finally everything boils down to the main estimate~:
\begin{theo}\cite{A, A2, AN2} \label{t:main}If the diameters of the supports of the functions $P_k$ are small enough, the following holds.

Choose $\chi$ smooth, compactly supported in $[E-\varepsilon, E+\varepsilon]$, and such that $\|\chi\|_\infty\leq 1$. For any $\tilde c>0$, there exists $\hbar_{\tilde c}>0$ such that, for all $\hbar<\hbar_{\tilde c}$, for
$n\leq \tilde c|\log\hbar|$, and
any pair of sequences $\alpha$, $\alpha'$ of length $n$,
\begin{equation}\label{e:stronglocali}
\left\|\pi_{\alpha'}(n)\pi_{\alpha} \,\chi(-\hbar^2\Delta) \right\|
 \leq C\hbar^{-(d-1)/2}\,e^{-n(d-1)\sqrt{E-\varepsilon}} .
\end{equation}
 (The constant $C$ is an absolute constant).
\end{theo}
We note that this result is an improvement of the estimate of \cite{A} (where the prefactor was only $\hbar^{-d/2}$) and \cite{AN2} (where the support of $\chi$ was assumed to shrink with $\hbar$). Proving Theorem \ref{t:entropy} using the weaker results of \cite{A} and \cite{AN2} turned out to be more painful than reproving Theorem \ref{t:main} directly. This proof is provided in \cite{A2}, section 5.
In what follows, the integer $n$ will always be taken equal to $\lfloor\tilde c|\log\hbar|\rfloor$, where $\tilde c$ will be fixed in the next section. We assume that $L$ is large enough so that $\hbar^{L-K_0}$
is negligible in comparison with $\hbar^{-(d-1)/2}\,e^{-n(d-1)\sqrt{E-\varepsilon}}$. Under all these assumptions, we have~:
\begin{prop}\label{p:WEUP}
Let $(\phi_\hbar)_{\hbar\to 0}$ be a sequence of normalized states satisfying the assumptions of proposition~\ref{c:WEUP} with $L$ large enough so that $\hbar^{L-K_0}$
is negligible in comparison with $\hbar^{-(d-1)/2}\,e^{-n(d-1)\sqrt{E-\varepsilon}}$ for $n=\lfloor\tilde c|\log\hbar|\rfloor$.
Then, in the semiclassical limit, the entropies of $\phi_\hbar$ at time $n=\lfloor\tilde c|\log\hbar|\rfloor$
satisfy
\begin{equation}\label{e:ineg}
\frac{h^+_n(U^{n\hbar}\phi_\hbar)+h_n^-( \phi_\hbar)}{2n}\geq (d-1)\sqrt{E-\varepsilon}- \frac{(d-1)}{2\tilde c}+\cO(n^{-1})
.
\end{equation}
\end{prop}

\subsection{Subadditivity until the Ehrenfest time}\label{s:subadd}
In this paragraph, we fix a sequence of normalized states $(\phi_\hbar)_{\hbar\to 0}$ satisfying \eqref{e:local} ($\chi$ is always assumed to be supported in $[E-\varepsilon, E+\varepsilon]$). We fix some arbitrary $\delta>0$, and introduce the Ehrenfest time,
\begin{equation}\label{e:Ehrenf}
n_{Ehr}(\hbar,E,\varepsilon):=\left\lfloor\frac{(1-\delta)|\log \hbar|}{\sqrt{E+\varepsilon }}\right\rfloor\,.
\end{equation}
\begin{rema}
The Ehrenfest time is the largest time interval on which the (non-commutative) dynamical system formed by the flow $(U^{t\hbar})$ acting on pseudodifferential operators (supported in $\{\|\xi\|^2\in[E-\varepsilon, E+\varepsilon]\}$) is commutative, up to small errors going to $0$ with $\hbar$.
\end{rema}

We take $n=n_{Ehr}(\hbar,E,\varepsilon)$ (in other words, we take $\tilde c=\frac{(1-\delta)}{\sqrt{E+\varepsilon }}$), and 
we use a subadditivity property of the entropies $h_n^+$ and $h_n^-$ to go from \eqref{e:ineg} for $n= n_{Ehr}(\hbar,E,\varepsilon)$ to a fixed, arbitrary, integer $n_0$. The proof of the next proposition is given in \cite{AN2}
in the case when $\phi_\hbar$ is an eigenfunction of $\Delta$. It can easily be adapted to the case of an arbitrary $\phi_\hbar$ and yields~:
\begin{prop} [Subadditivity]\label{p:subadd} Let $E\geq 0$ and $\varepsilon>0$.
For $\delta>0$ arbitrary, define the Ehrenfest time $n_{Ehr}(\hbar, E, \varepsilon)$ as in \eqref{e:Ehrenf}.
Let $(\phi_\hbar)_{\hbar\rightarrow0}$ be a normalized family satisfying \eqref{e:local}, where $\chi$ is supported in $[E-\varepsilon, E+\varepsilon]$, and $L$ is chosen large enough.

 For any $n_0\in\mathbb{N}$, there exists a positive $R_{n_0}(\hbar)$, with $R_{n_0}(\hbar)\rightarrow0$ as $\hbar\rightarrow0$, such that
for any $\hbar\in(0,1]$, any
$n_0,m\in\mathbb{N}$ with $n_0+m\leq n_{Ehr}(\hbar)$, we have
$$
h^+_{n_0+m}(\phi_\hbar)\leq h^+_{m}(\phi_\hbar)+h^+_{n_0}(U^{m\hbar}\phi_\hbar)
+R_{n_0}( \hbar)\, ,
$$
$$
h^-_{n_0+m}(\phi_\hbar)\leq h^-_{n_0}(\phi_\hbar)+h^-_{m}(U^{n_0\hbar}\phi_\hbar)
+R_{n_0}( \hbar).
$$
\end{prop}
Let $n_0\in\mathbb{N}$ be fixed and $n=n_{Ehr}(\hbar,E,\varepsilon)$.
Using the Euclidean division $n=qn_0+r $, with $r< n_0$, Proposition~\ref{p:subadd} implies that
for $\hbar$ small enough,
$$
\frac{h_n^+(\phi_\hbar)}{n}\leq \frac{\sum_{k=0}^{q-1} h_{n_0}^+(U^{kn_0\hbar}\phi_\hbar)}{qn_0}
+\frac{h_r^+(U^{qn_0\hbar} \phi_\hbar)}{n}+\frac{R_{n_0}(\hbar)}{n_0}
$$
and
$$
\frac{h_n^-(\phi_\hbar)}{n}\leq \frac{\sum_{k=0}^{q-1} h_{n_0}^-(U^{(r+kn_0)\hbar}
 \phi_\hbar)}{qn_0}
+\frac{h_r^-(U^{r\hbar}\phi_\hbar)}{n}+\frac{R_{n_0}(\hbar)}{n_0}.
$$
Note that $h_r^+(U^{qn_0\hbar} \phi_\hbar)+h_r^-(U^{r\hbar}\phi_\hbar)$
stays uniformly bounded (by $\log n_0$) when $\hbar\to 0$.  Combining the subadditivity property with Proposition \ref{p:WEUP}, we find that
\begin{equation}\label{e:fixedtime}
\frac{\sum_{k=0}^{q-1} \left(  h_{n_0}^+(U^{kn_0\hbar}U^{n\hbar}\phi_\hbar)+  h_{n_0}^-(U^{(r+kn_0)\hbar}
 \phi_\hbar)\right)}{2qn_0}
\geq (d-1)\sqrt{E-\varepsilon}- \frac{(d-1)\sqrt{E+\varepsilon}}{2(1-\delta)}-\frac{R_{n_0}( \hbar)}{n_0}+\cO_{n_0}(1/n),
\end{equation}
for $n=n_{Ehr}(\hbar,E,\varepsilon)$.

\subsection{The conclusion}
The interval $[E_1, E_2]$ is fixed. Consider $E$ in $[E_1,E_2]$ and a sequence of normalized states $(u_{\hbar})_{\hbar\rightarrow0}$ that satisfies
\eqref{e:local-etats}. We may assume without loss of generality that $\bbbone_{[E_1,E_2]}\left(-\hbar^2\Delta\right) u_{\hbar}=u_\hbar$ (since the semiclassical limits associated with $u_\hbar$ and $\bbbone_{[E_1,E_2]}\left(-\hbar^2\Delta\right) u_{\hbar}$ will be the same). We fix a function $\chi\in \cC_o^\infty(\R)$, supported in $[-1, 1]$ such that $\sum_{k\in\Z}\chi^2(x-k)\equiv 1$. For $N\in\N$, we denote $\varepsilon=\frac{E_2-E_1}N$, and $\chi_{j} (x)=\chi\left(\frac{x-E_1-j\varepsilon}\varepsilon\right)$ ($j=0,\ldots, N$). We have $u_\hbar=\sum_{j=0}^N\chi^2_j\left(-\hbar^2\Delta\right) u_{\hbar}$ and thus $\|u_\hbar\|^2=\sum_{j=0}^N\|\chi_j\left(-\hbar^2\Delta\right) u_{\hbar}\|^2$. We will denote $u_j=\chi_j\left(-\hbar^2\Delta\right) u_{\hbar}$ and $\tilde u_j=\frac{u_j}{\|u_j\|}$.
For $t\in\R$, we apply \eqref{e:fixedtime} to $\phi_\hbar=U^t\tilde u_j$ and obtain
\begin{multline}\label{e:fixedtimeenergy}
\frac{\sum_{k=0}^{q-1} \left(  h_{n_0}^+(U^{kn_0\hbar}U^{n\hbar}U^t\tilde u_j)+  h_{n_0}^-(U^{(r+kn_0)\hbar}U^t\tilde u_j)\right)}{2qn_0}\\
\geq (d-1)\sqrt{E_1+(j-1)\varepsilon}- \frac{(d-1)}{2(1-\delta)}\sqrt{E_1+(j+1)\varepsilon}-\frac{R_{n_0}( \hbar)}{n_0}+\mathcal{O}_{n_0}(1/|\log\hbar|),
\end{multline}
If we multiply by $\theta(t)$ (satisfying $\theta\in L^1(\R,\R_+)$ and $\int\theta=1$), integrate with respect to $t$, and take into account the fact that $(kn_0+r)\hbar\To 0$ and $n\hbar\To 0$, we find that
\begin{multline}\label{e:fixedtimes-locstates}\int\theta(t)\frac{h_{n_0}^+(U^t\tilde u_j)+
 h_{n_0}^-(U^t\tilde u_j)}{2n_0}dt \\ \geq (d-1)\sqrt{E_1+(j-1)\varepsilon}- \frac{(d-1)}{2(1-\delta)}\sqrt{E_1+(j+1)\varepsilon}+o_{n_0,N}(1) .\end{multline}

We now use the fact that
$$\sum_{j=0}^N\left\|\pi_{\alpha}U^t u_j\right\|^2=\left\|\pi_{\alpha}U^t u_{\hbar}\right\|^2+o_{n_0,N}(1),$$
where the remainder is uniform for $\alpha$ of length $n_0$, and for $t\in\R$. This comes from the fact that $\pi_{\alpha}$ and $\chi_j\left(-\hbar^2\Delta\right)$ both belong to $\Psi^{0, 0} (M)$, from which we have the commutator estimate $\|[ \pi_{\alpha}, \chi_j\left(-\hbar^2\Delta\right)] \|=\cO(\hbar)$ (we also use the fact that $\sum \chi_j^2=1$).
 
If we combine this with the concavity of the map $x\mapsto-x\log x$, we can verify that as $\hbar$ tends to $0$,
$$ \sum_{j=0}^N\left\| u_j\right\|^2h_{n_0}^+\left(U^t\tilde u_j\right)\leq h_{n_0}^+\left(U^tu_{\hbar} \right)+o_{n_0,N}(1),$$
where the remainder is uniform for $t$ in $\R$. A similar inequality holds with $h_{n_0}^-$. Combined with  \eqref{e:fixedtimes-locstates}, this yields that
\begin{multline}\label{e:before}\int\theta(t)\frac{h_{n_0}^+(U^tu_\hbar)+
 h_{n_0}^-(U^t u_\hbar)}{2n_0}dt \\ \geq  \sum_{j=0}^N
 \left\| u_j\right\|^2\left[ (d-1)\sqrt{E_1+(j-1)\varepsilon}- \frac{(d-1)}{2(1-\delta)}\sqrt{E_1+(j+1)\varepsilon}\right]+
 o_{n_0}(1) .\end{multline}
 We define the following averaged entropy
  \begin{align}\label{e:averentropy}
h^-_n(\phi,\theta)&=-\sum_{|\alpha|=n}\left(\int\theta(t)\left\|\pi^*_{\alpha}U^t \phi\right\|^2dt \right)\log\left(\int\theta(t)\left\|\pi^*_{\alpha}U^t \phi\right\|^2dt \right),\\
h^+_n(\phi,\theta)&=-\sum_{|\alpha|=n}\left(\int\theta(t)\left\|\pi_{\alpha}U^t \phi\right\|^2dt \right)\log\left(\int\theta(t)\left\|\pi_{\alpha}U^t \phi\right\|^2dt \right).
\end{align}

Using again the concavity of $x\mapsto -x\log x$, \eqref{e:before} implies
\begin{multline}\label{e:after} \frac{h_{n_0}^+(u_\hbar, \theta)+
 h_{n_0}^-(u_\hbar, \theta)}{2n_0} \\ \geq  \sum_{j=0}^N
 \left\| u_j\right\|^2\left[ (d-1)\sqrt{E_1+(j-1)\varepsilon}- \frac{(d-1)}{2(1-\delta)}\sqrt{E_1+(j+1)\varepsilon}\right]+o_{n_0,N}(1) .\end{multline}

 We can now take the limit $\hbar\To 0$. If the semiclassical measure associated with the family
 $(U^tu_\hbar)$ decomposes as $\mu_t=\int \mu_{t, E} d\nu(E)$, then $\|u_j\|^2$ converges to $\int \chi_j^2(E) d\nu(E)$. On the left-hand side of \eqref{e:after}, $h_{n_0}^+(u_\hbar, \theta)$ and
$ h_{n_0}^-(u_\hbar, \theta)$ both converge to
$$-\sum_{|\alpha|=n}\left(\int\theta(t) \mu_t((P^2_{\alpha_{n-1}}\circ g^{n-1})\ldots
(P^2_{\alpha_{1}}\circ g^{1})P^2_{\alpha_{0}}
)dt \right)\log\left(\int\theta(t) \mu_t((P^2_{\alpha_{n-1}}\circ g^{n-1})\ldots
(P^2_{\alpha_{1}}\circ g^{1})P^2_{\alpha_{0}}
)dt  \right).$$
 
 After taking the limit $\hbar\To 0$,
 we take the limit $N\To +\infty$~: on the right-hand side of \eqref{e:after}, this transforms the discrete lower bound $$\sum_{j=0}^N
   \left[(d-1)\sqrt{E_1+(j-1)\varepsilon}- \frac{(d-1)}{2(1-\delta)}\sqrt{E_1+(j+1)\varepsilon}\right]\int \chi_j^2(E) d\nu(E)$$ into the integral $\frac{d-1}{2}\int \sqrt{E} d\nu(E)$ (after taking $\delta\rightarrow0$).

   Finally we let $n_0\To+\infty$, which allows to go from $h^\pm_{n_0}$ on the left-hand side of \eqref{e:after} to the Kolmogorov-Sinai entropy $h_{KS}$; for this step, details can be found in \cite{AN2} (paragraph $2.2.8$).

 At this stage, we obtain
$$h_{KS}\left( \int \theta(t) \mu_tdt\right)=h_{KS}\left( \int \theta(t) \mu_{t, E} d\nu(E)dt\right)\geq\frac{d-1}{2}\int_{E_1}^{E_2} \sqrt{E}{d\nu(E)} .$$

If we use the same argument, replacing $u_\hbar$ by $f(-\hbar^2\Delta)u_\hbar$ (where $f$ is a smooth function on $[E_1,E_2]$ such that $\int f^2(E)d\nu(E)=1$), we obtain by the same argument
$$h_{KS}\left( \int \theta(t) \mu_{t, E}f^2(E) d\nu(E)dt\right)\geq\frac{d-1}{2}\int_{E_1}^{E_2} \sqrt{E}f^2(E){d\nu(E)} .$$

We finally use the fact that $\mu\mapsto h_{KS}(\mu)$ is affine, to convert the last inequality into
$$\int_{\R}\int_{E_1}^{E_2}\theta(t)f^2(E)\left(h_{KS}(\mu_{t,E})-\frac{d-1}{2}\sqrt{E}\right)dtd\nu(E)\geq 0;$$
this inequality holds for all $\theta$ in $L^1(\mathbb{R},\mathbb{R}_+)$ such that $\int\theta=1$ and $f$ in $\mathcal{C}^{\infty}_o(\mathbb{R}_+,\mathbb{R})$
 such that  $\int f^2(E)d\nu(E)=1$. As a consequence, one has for $Leb\otimes \nu$-almost every $(t, E)$, 
$$h_{KS}(\mu_{t,E})\geq\frac{d-1}{2}\sqrt{E}.\square$$

\appendix

\section{From entropy estimates to observability}\label{s:control}

In this short appendix, we explain how we can go from the entropy estimates of Theorem~\ref{t:entropy} to the observability estimate of Theorem \ref{t:observability}. According to Lebeau~\cite{Leb}, it is sufficient to prove the following weak observability result to deduce Theorem~\ref{t:observability}:
\begin{theo}\label{t:weakobservability} 
Under the assumptions of Theorem~\ref{t:observability}, for all $T>0$, there exists $C_{T, a}>0$ such that, for all $u$~:
\begin{equation}\label{e:wobservability}\|u\|^2_{L^2(M)}\leq C_{T, a}\left(\int_0^T\| ae^{\imath t\frac\Delta2}u\|^2_{L^2(M)} dt+\|u\|_{H^{-1}(M)}^2\right).\end{equation}
\end{theo}
For the sake of completeness, we briefly recall the argument of Lebeau to deduce observability from a weak observability estimate at time $T$. First, for $T'>T$, we introduce the subspace
$$N(T'):=\left\{\varphi\in L^2(M):\forall 0\leq t\leq T',\ a(x)(e^{\imath t\Delta}\varphi)(x)=0\right\}.$$
From weak observability and the compactness of the injection $L^2\subset H^{-1}$, we can deduce that for $T'>T$, this subspace is finite dimensional. One can also verify that for every $T<T''<T'$ and for every $\varphi$ in $N(T')$, $\Delta\varphi$ belongs to $N(T'')$ (by taking the limit of the sequence $\frac{e^{\imath\epsilon\Delta}\varphi-\varphi}{\epsilon}$, which belongs to $N(T'')$ for $\epsilon$ small enough, and is bounded in $H^{-2}(M)$). 

This implies that $\Delta$ is an operator from the finite dimensional subspace $N(T')$ into itself. As $a$ is nontrivial, one can deduce the existence of an eigenfunction of the laplacian which is equal to $0$ on a nonempty open set. By Aronszajn-Cordes' theorem~\cite{Ho} (section $17.2$), this eigenfunction is necessarly $0$ and the subspace $N(T')$ is reduced to $\{0\}$. By contradiction, we can finally deduce that observability holds for $T'>T$.\\

In order to prove Theorem~\ref{t:weakobservability}, we proceed by contradiction and make the assumption that there exist a sequence of normalized vectors $(u_{n})_{n\in\mathbb{N}}$ in $L^2(M)$ and $T>0$ such that
\begin{equation}\label{e:contradiction}\lim_{n\rightarrow+\infty}\left(\int_0^T\| ae^{\imath t\frac\Delta2}u_{n}\|^2_{L^2(M)} dt+\|u_n\|_{H^{-1}(M)}^2\right)=0.
\end{equation}
This implies that $u_n$ converges to $0$, weakly in $L^2$.
For every $t$ in $\mathbb{R}$, we introduce the distribution
$$  \mu_n(t)(b):=\langle u_n|e^{-\imath t\frac\Delta2}\Op_1(b)e^{\imath t\frac\Delta2}u_n\rangle_{L^2(M)},$$
defined for all $b\in\cS_0$ (see the remark following Theorem~\ref{t:entropy}).
This distribution is the analogue of~(\ref{general-def-distrib}) in the microlocal setting~\cite{Ge}. As before, the map $t\mapsto \mu_n(t)$ belongs to $L^{\infty}(\mathbb{R},\cS_0')$. Thus, there exists a subsequence $(u_{n_k})_k$ and $\mu$ in $L^{\infty}(\mathbb{R}, \cS_0')$ such that
\[
\int_{\R\times \widehat{T^*M}}%
\theta(t)b(x,\xi)\mu_{n_k}(t)(dx,d\xi)dt\Lim_{k\To+\infty}\int_{\R\times \widehat{T^*M}}%
\theta(t)b(x,\xi)\mu(t)(dx,d\xi)dt
\]
for all $\theta\in L^{1}(\R)$ and $b\in  \cS_0.$ Besides, as above, $\mu(t)$ is a probability measure on the compactified cotangent bundle $\widehat{T^*M}$, and is invariant under the normalized geodesic flow.
As $u_n(t)=e^{\imath t\frac\Delta2}u_n$ converges weakly to $0$ for every $t$ in $\mathbb{R}$,
each $\mu(t)$ is actually supported at infinity, and may thus be identified with a probability measure
on the unit sphere bundle $S^*M$, invariant under the geodesic flow.

From Theorem~\ref{t:entropy} and the associated remark, we know that for almost every $t$ in $\mathbb{R}$, $h_{KS}(\mu(t))\geq \frac{d-1}{2}$. We will now use the fact that the Hausdorff dimension of $K_a$ is less than $d$. From~\cite{PeSa} (Theorem $4.2$), this implies that the topological entropy satisfies
$$h_{top}(K_a,(g^t)):=\sup_{\mu\in\mathcal{M}(S^*M,g^t)}\left\{h_{KS}(\mu):\mu(K_a)=1\right\}<\frac{d-1}{2}.$$
Using property~(\ref{e:contradiction}), we know that $\int_{S^*M\times[0,T]}a^2(x,\xi)\mu(t)(dx,d\xi)dt=0$. In particular, this implies that $\mu(t)(S^*M\backslash K_a)=0$ for almost every $t$ in $[0,T]$ (as $\mu(t)$ is $g^s$-invariant) and it leads to a contradiction.$\square$

\section{Pseudo-differential calculus on a manifold}\label{pdo-manifold}

In this section, we recall some facts of pseudodifferential calculus; details can be found in~\cite{EZ}.
We define on $\mathbb{R}^{2d}$ the following class of (semiclassical) symbols:
\begin{multline*}S^{m,k}(\mathbb{R}^{2d}):=\{a=a_\hbar\in C^{\infty}(\mathbb{R}^{2d}): \forall K\subset \mathbb{R}^{d}\; {\rm compact },\\ \forall \alpha, \beta, \exists C_{\alpha,\beta},\forall (x, \xi)\in K\times\R^d, |\partial^{\alpha}_x\partial^{\beta}_{\xi}a|\leq C_{\alpha,\beta}\hbar^{-k}\langle\xi\rangle^{m-|\beta|}\}.
\end{multline*}
Let $M$ be a smooth compact riemannian $d$-manifold without boundary. Consider a finite smooth
atlas $(f_l,V_l)$ of $M$, where each $f_l$ is a smooth diffeomorphism from the open subset $V_l\subset M$ to a bounded open set $W_l\subset\mathbb{R}^{d}$. To each $f_l$ correspond a pull back $f_l^*:C^{\infty}(W_l)\rightarrow C^{\infty}(V_l)$ and a canonical map $\tilde{f}_l$ from $T^*V_l$ to $T^*W_l$:
$$\tilde{f}_l:(x,\xi)\mapsto\left(f_l(x),(Df_l(x)^{-1})^T\xi\right).$$
Consider now a smooth locally finite partition of identity $(\phi_l)$ adapted to the previous atlas $(f_l,V_l)$. That means $\sum_l\phi_l=1$ and $\phi_l\in C_o^{\infty}(V_l)$. Then, any observable $a$ in $C^{\infty}(T^*M)$ can be decomposed as follows: $a=\sum_l a_l$, where $a_l=a\phi_l$. Each $a_l$ belongs to $C^{\infty}(T^*V_l)$ and can be pushed to a function $\tilde{a}_l=(\tilde{f}_l^{-1})^*a_l\in C^{\infty}(T^*W_l)$. As in~\cite{EZ}, define the class of symbols of order $m$ and index $k$:
\begin{equation}\label{defpdo}S^{m,k}(T^{*}M):= \left\{a=a_\hbar\in C^{\infty}(T^*M):  \forall \alpha, \beta, \exists C_{\alpha,\beta},|\partial^{\alpha}_x\partial^{\beta}_{\xi}a|\leq C_{\alpha,\beta}\hbar^{-k}\langle\xi\rangle^{m-|\beta|}\right\}.\end{equation}
Then, for $a\in S^{m,k}(T^{*}M)$ and for each $l$, one can associate to the symbol $\tilde{a}_l\in S^{m,k}(\mathbb{R}^{2d})$ the standard Weyl quantization:
$$\Op_{\hbar}^{w}(\tilde{a}_l)u(x):=\frac{1}{(2\pi\hbar)^d}\int_{\mathbb{R}^{2d}}e^{\frac{\imath}{\hbar}\langle x-y,\xi\rangle}\tilde{a}_l\left(\frac{x+y}{2},\xi;\hbar\right)u(y)dyd\xi,$$
where $u\in\cC_o^\infty(\mathbb{R}^d)$. Consider now a smooth cutoff $\psi_l\in C_c^{\infty}(V_l)$ such that $\psi_l=1$ close to the support of $\phi_l$. A quantization of $a\in S^{m,k}(T^*M)$ is then defined in the following way:
\begin{equation}\label{pdomanifold}\Op_{\hbar}(a)(u):=\sum_l \psi_l\times\left(f_l^*\Op_{\hbar}^w(\tilde{a}_l)(f_l^{-1})^*\right)\left(\psi_l\times u\right),\end{equation}
where $u\in C^{\infty}(M)$. According to the appendix of~\cite{EZ}, the quantization procedure $\Op_{\hbar}$ sends $S^{m,k}(T^{*}M)$ onto the space of pseudodifferential operators of order $m$ and of index $k$, denoted $\Psi^{m,k}(M)$. It can be shown that the dependence in the cutoffs $\phi_l$ and $\psi_l$ only appears at order $2$ in $\hbar$ and the principal symbol map $\sigma_0:\Psi^{m,k}(M)\rightarrow S^{m,k}/S^{m,k-1}(T^{*}M)$ is then intrinsically defined. 

At various places in this paper, a larger class of symbols should be considered, as in~\cite{DS} or~\cite{EZ}. For $0\leq\nu<1/2$:
\begin{equation}\label{symbol}S^{m,k}_{\nu}(T^{*}M)= \left\{a=a_\hbar\in C^{\infty}(T^*M):\forall \alpha, \beta, \exists C_{\alpha,\beta}, |\partial^{\alpha}_x\partial^{\beta}_{\xi}a|\leq C_{\alpha,\beta}\hbar^{-k-\nu|\alpha+\beta|}\langle\xi\rangle^{m-|\beta|} \right\}.\end{equation}
Results of~\cite{DS} can be applied to this new class of symbols. For example, if $M$ is compact, a symbol of $S^{0,0}_{\nu}$ gives a bounded operator on $L^2(M)$ (with norm independent of $\hbar\leq 1$).

Even if the Weyl procedure is a natural choice to quantize an observable $a$ on $\mathbb{R}^{2d}$, it is sometimes preferrable to use a quantization that satisfies the additional property~: $\Op_{\hbar}(a)\geq 0$ if $a\geq0$. This can be achieved thanks to the anti-Wick procedure, see~\cite{HeMaRo}.
 For $a$ in $S^{0,0}_{\nu}(\mathbb{R}^{2d})$, that coincides with a function on $M$ outside a compact subset of $T^*M$, one has
\begin{equation}\label{equivalence-positive-quantization}\|\Op_{\hbar}^w(a)-\Op_{\hbar}^{AW}(a)\|_{L^2}\leq C\sum_{|\alpha|\leq D}\hbar^{\frac{|\alpha|+1}{2}}\|\partial^{\alpha}da\|,
\end{equation}
where $C$ and $D$ are some positive constants that depend only on the dimension $d$.
To get a positive procedure of quantization on a manifold, one can replace the Weyl quantization by the anti-Wick one in definition~(\ref{pdomanifold}). We will denote $\Op_{\hbar}^+(a)$ this new choice of quantization, well defined for every element in $S^{0,0}_{\nu}(T^*M)$ of the form $b(x)+c(x,\xi)$ where $b$ belongs to $S^{0,0}_{\nu}(T^*M)$ and $c$ belongs to $\mathcal{C}^{\infty}_o(T^*M)\cap S^{0,0}_{\nu}(T^*M)$. We underline the fact that $\Op_{\hbar}^+(1)=\text{Id}_{L^2(M)}$.

\end{document}